\documentclass[letterpaper,10pt]{article}
\usepackage[textwidth=11.5cm, textheight=22cm, top=2.0cm, centering]{geometry}
\usepackage{setspace}
\usepackage[T1]{fontenc}
\usepackage{lmodern}
\usepackage[english]{babel}
\usepackage{microtype} 
\usepackage[scaled=0.95]{inconsolata} 
\makeatletter
\renewcommand{\texttt}[1]{{\ttfamily #1}}

\renewcommand{\mathtt}[1]{\text{\texttt{#1}}} 
\makeatother
\usepackage{csquotes}
\usepackage{etoolbox}
\usepackage{bbding}
\usepackage{changepage}
\usepackage{booktabs}
\usepackage{graphicx}
\usepackage[x11names, dvipsnames]{xcolor}
\definecolor{Linkz}{RGB}{30, 110, 170}
\definecolor{Darkenta}{RGB}{185, 35, 90}
\definecolor{Magentz}{RGB}{255, 35, 170}
\definecolor{Lightenta}{RGB}{254, 232, 255}
\definecolor{Reference}{RGB}{35, 180, 90}
\definecolor{Periwinkle}{RGB}{102, 51, 255}
\definecolor{yello}{RGB}{255, 245, 230}
\definecolor{Greeno}{RGB}{0, 140, 100}
\definecolor{Leeno}{RGB}{239, 255, 232}
\definecolor{Nicegreen}{RGB}{100, 200, 130}
\usepackage{amsmath, amsthm, amscd, amssymb, stmaryrd}
\usepackage{stackengine}
\usepackage{bussproofs}
\usepackage[square]{natbib}

\usepackage[
colorlinks=true,
linkcolor=Darkenta,
citecolor=Greeno,
urlcolor=Darkenta,
]{hyperref}
\usepackage{titlesec} 
\usepackage{abstract} 
\usepackage{enumitem} 
\usepackage{ccicons} 
\usepackage{float}
\usepackage{placeins} 
\usepackage{tikz} 
\usepackage{tikz-cd}
\usetikzlibrary{positioning,fit,patterns}
\usetikzlibrary{arrows.meta} 
\usepackage{adjustbox} 
\usepackage{pgfmath} 
\usepackage{ifthen} 

\newtheoremstyle{upright}
{6pt plus 2pt minus 2pt} 
{6pt plus 2pt minus 2pt} 
{\normalfont} 
{} 
{\bfseries} 
{.} 
{.5em} 
{} 
\theoremstyle{upright}
\theoremstyle{upright}
\newtheorem{theorem}{Theorem}[subsection]
\newtheorem{remark}[theorem]{Remark}

\newtheorem{definition}[theorem]{Definition}

\newtheorem{proposition}[theorem]{Proposition}

\makeatletter
\renewenvironment{proof}[1][Proof]{%
	\par\pushQED{\qed}%
	\normalfont
	\topsep6\p@\@plus6\p@\relax
	\trivlist
	\item[\hskip\labelsep\slshape #1\@addpunct{.}]%
}{%
	\popQED\endtrivlist\@endpefalse
}

\makeatother

\usepackage{algpseudocode}
\usepackage[most]{tcolorbox}
\usepackage{listings}
\newtcolorbox{breakbox}[2][]{%
	breakable,
	={#2},
	fonttitle=\bfseries,
	colback=white,
	colframe=black!20,
	coltitle=black,
	colbacktitle=white,
	boxrule=0.5pt,
	arc=0pt,
	boxsep=7pt,
	left=3pt,
	right=2pt,
	top=2pt,
	bottom=4pt,
	fontupper=\small\sffamily, 
	#1
}

\AtBeginEnvironment{quotation}{\sffamily}

\newcommand{\customsectionstyle}[2]{%
	\titleformat{\section}[block]
	{\normalfont\fontsize{#1}{1.2\dimexpr#1\relax}\selectfont\centering}
	{\thesection}{1em}%
	{%
		\ifthenelse{\equal{#2}{true}}{\MakeUppercase}{\relax}%
	}%
}
\newcommand{\customsectionspacing}[3]{%
	\titlespacing*{\section}{#1}{#2}{#3}%
}
\newcommand{\customsubsectionstyle}[2]{%
	\titleformat{\subsection}[block]
	{\normalfont\fontsize{#1}{1.2\dimexpr#1\relax}\selectfont\centering}
	{\thesubsection}{1em}%
	{%
		\ifthenelse{\equal{#2}{true}}{\MakeUppercase}{\relax}%
	}%
}
\newcommand{\customsubsectionspacing}[3]{%
	\titlespacing*{\subsection}{#1}{#2}{#3}%
}
\customsectionstyle{14pt}{true}
\customsectionspacing{0pt}{3.5ex plus 1ex minus 0.2ex}{2.5ex}
\customsubsectionstyle{11pt}{true}
\customsubsectionspacing{0pt}{4ex plus 1ex minus 0.2ex}{2ex}
\usepackage{caption}

\makeatletter
\newcommand{\shorttitle}[1]{\def\@shorttitle{#1}}
\newcommand{\email}[1]{\def\@email{#1}}
\newcommand{\metadata}[1]{\def\@metadata{#1}}
\renewcommand{\maketitle}{%
	\begin{center}
		\vfill
		{\fontsize{18pt}{19pt}\selectfont \@title \par}
		\vspace{1em}
		{\normalsize \@author \par}
		\vspace{0.1em}
		{\normalsize \@date \par}
	\end{center}
}
\makeatother

\begin{document}

\title{\uppercase{ON THE GOLDEN RATIO\\AND STABLE SELF-APPLICATION}}

\author{Milan Rosko}
\date{June 2026}

\maketitle

\begin{center}\footnotesize{
		ORCID: \href{https://orcid.org/0009-0003-1363-7158}{\footnotesize\textsf{0009-0003-1363-7158}}\\
}
\end{center}

\begin{abstract} \vspace{-1ex}\footnotesize{
This paper studies a boundary between local self-application and global self-certification. Irrational quantities are treated operationally, as procedures whose approximations are refined by effective update rules. The golden ratio $\Phi$ is used as a model of stable local recurrence: the reciprocal update $R(x)=1+1/x$ has a unique positive fixed point and admits finite witnessed approximations. By contrast, global reflection asks a system to certify its own correctness uniformly. The proof-theoretic claim is therefore contrastive: primitive-recursive proof checking and local soundness preserve correctness through bounded checks and bounded witnesses, but they do not yield internal global reflection. No complexity advantage, decision procedure, or new reflection principle is claimed.

	}
\end{abstract}

\begin{figure}[H]
	\centering
	\includegraphics[width=0.78\textwidth]{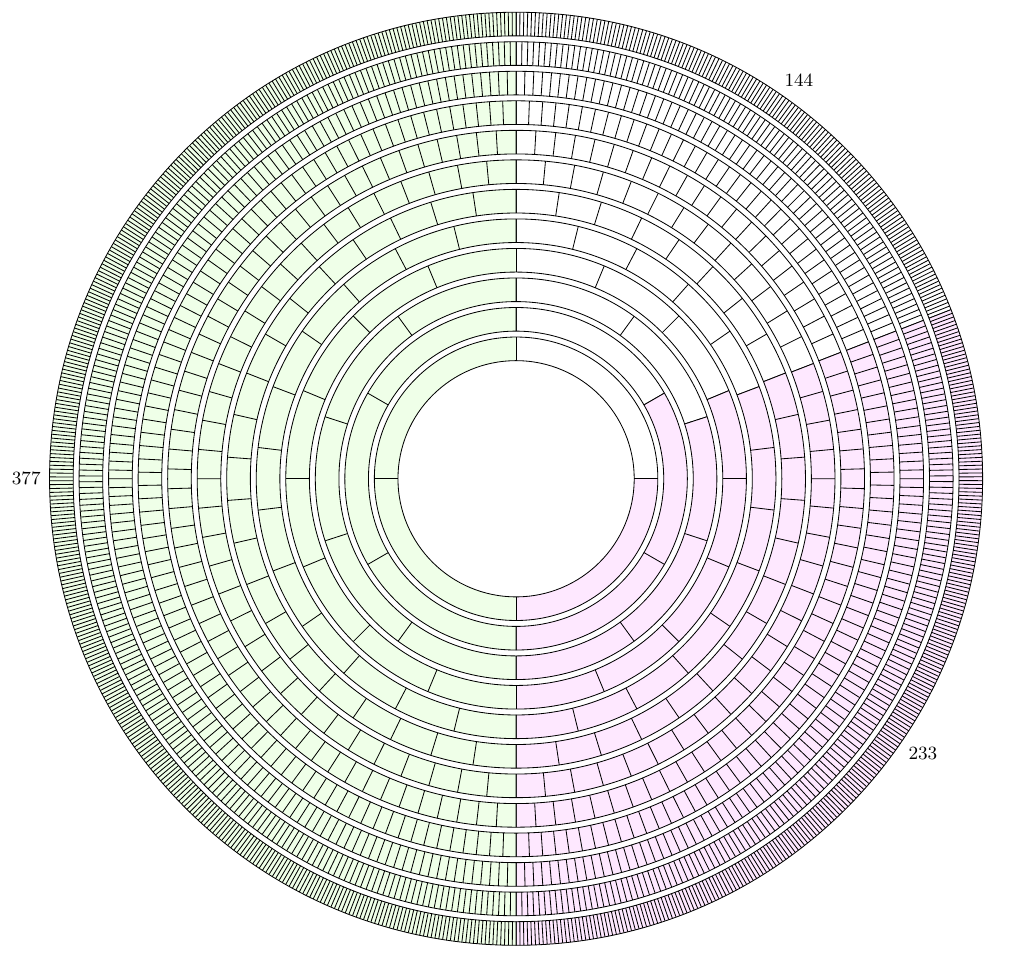}
	\caption{Geometric intuition for Fibonacci alignment. Each annulus is partitioned by a consecutive triple $(F_n,F_{n+1},F_{n+2})$. Since $F_n+F_{n+1}=F_{n+2}$, the $F_{n+2}$ sector occupies one half of the annulus, while the opposite half is split into adjacent sectors of sizes $F_n$ and $F_{n+1}$. Across the nested annuli this repeats the same Fibonacci pattern at successive scales. The $\Phi$ approximation is read from the adjacent terms $F_n$ and $F_{n+1}$.}
	\label{fig:wheel}
\end{figure}

\section{Introduction}\label{sec:introduction}

\subsection{The extreme and mean ratio}\label{subsec:extreme-mean-ratio}

Programmatic papers have a specific obligation: they must distinguish, with care, between what they establish and what they propose. The difficulty in the present case is that self-application, slow convergence, and inferential systems are not unrelated themes: there is a genuine structural resonance among them. Yet this resonance is liminal. It is strong enough to organize a useful boundary model, but not strong enough to license a new proof system, a complexity claim, or an arithmetic reduction of logic.

The golden ratio, moreover, is a \emph{difficult} constant to invoke, especially in the context of proof theory and weak systems. On the one hand, it has often been surrounded by exaggerated claims, particularly where aesthetic, biological, or cosmological patterns are inferred from approximate numerical resemblance, leaving a tableau of misconceptions as a warning, as entertainingly shown by \citet{markowsky92}. On the other hand, it does not follow that $\Phi$ is mathematically inert, nor that natural accounts of it are spurious as a general rule.

\subsection{Roadmap and contribution}\label{subsec:roadmap}

The contribution is a deliberately limited boundary model. Section~\ref{sec:definitions} defines resonance as an analogy between local application and stable proportional self-similarity, then separates the number-theoretic role of $\Phi$ from the finite coding role of Fibonacci recurrence. Section~\ref{sec:proof-checking} gives the local proof-checking substrate: line codes, axiom-head recognition, \textsc{MP} alignment, and bounded witnesses. Section~\ref{sec:verification-scope} states what this substrate establishes: primitive-recursive checking and local soundness, but not internal global reflection. Section~\ref{sec:discussion} collects the philosophical point: $\Phi$ models disciplined local self-application, while reflection marks the point where self-application is no longer justified by local witnessed checking alone.

\section{Definitions}\label{sec:definitions}

\subsection{Self-Similarity}\label{subsec:self-similarity}

\begin{definition}[Proportional resonance]
	The logical rule \textsc{Modus Ponens} has the local application form
	\begin{equation}
		A,\quad A\to B \quad \Longrightarrow \quad B
	\end{equation}
	in which the independent premise $A$ matches the antecedent occurrence inside $A\to B$, and detachment yields $B$. The proportional analogue is defined as follows. Let $a,b>0$ be magnitudes, let $\rho>0$, and define
	\begin{equation}
		m_{\rho}(x)=\frac{x}{\rho}.
	\end{equation}
	Write
	\begin{equation}
		a\xrightarrow{\rho} b
	\end{equation}
	when
	\begin{equation}
		m_{\rho}(a)=b.
	\end{equation}
	Equivalently,
	\begin{equation}
		a\xrightarrow{\rho} b \quad \Longleftrightarrow \quad \frac{a}{\rho}=b\; \Longleftrightarrow \quad \frac{a}{b}=\rho.
	\end{equation}
	Thus
	\begin{equation}
		a,\quad a\xrightarrow{\rho} b \quad \Longrightarrow \quad b
	\end{equation}
	is the proportional analogue of local application. A proportional resonance occurs when the same proportional morphism persists across adjacent scales:
	\begin{equation}
		a+b \xrightarrow{\rho} a \xrightarrow{\rho} b.
	\end{equation}
\end{definition}

\begin{proposition}[Golden resonance]
	If $a,b>0$ and
	\begin{equation}
		a+b \xrightarrow{\rho} a \xrightarrow{\rho} b,
	\end{equation}
	then the unique positive coherent ratio is
	\begin{equation}
		\rho=\Phi=\frac{1+\sqrt5}{2}.
	\end{equation}
\end{proposition}

\begin{proof}
	This means
	\begin{equation}
		\frac{a+b}{a}=\frac{a}{b}=\rho.
	\end{equation}
	Since $\rho=a/b$, we obtain
	\begin{equation}
		\rho = \frac{a+b}{a} = 1+\frac{b}{a} = 1+\frac{1}{\rho}.
	\end{equation}
	Hence
	\begin{equation}
		\rho^2=\rho+1.
	\end{equation}
	The unique positive solution is
	\begin{equation}
		\rho=\Phi=\frac{1+\sqrt5}{2}.
	\end{equation}
\end{proof}

\begin{remark}
	The proportional construction asks what happens when an analogue of this local application pattern is required to persist across adjacent scales. The answer is forced:
	\begin{equation}
		a+b \xrightarrow{\Phi} a \xrightarrow{\Phi} b.
	\end{equation}
	In the logical rule, the independent premise $A$ matches the antecedent occurrence inside $A\to B$, and detachment yields $B$. In the proportional analogue, an analogous local matching is required to persist across adjacent scales. This gives a model of disciplined local self-application, not a claim of equivalence and not a global reflection principle.
\end{remark}

\begin{remark}
	From now on, the role of $\Phi$ becomes narrow, operational and contrastive. It is not invoked as a hidden law of logic, nor as a numerological invariant of inference. $\Phi$ is the stable state of the reciprocal update
	\begin{equation}
		R(x)=1+\frac{1}{x},
	\end{equation}
	and this update repeats the same constructor at every finite stage of its continued fraction expansion. The point is not that logic is Fibonacci-shaped, but that $\Phi$ provides a precise example of local self-application whose stability can be witnessed from step to step.
\end{remark}

\begin{remark}
	Two arithmetical facts will therefore be kept separate. First, Fibonacci recurrence gives finite combinatorial structure, including canonical decompositions such as \textsc{Zeckendorf} representation. Second, the reciprocal update has $\Phi$ as its unique positive fixed point and generates Fibonacci ratios under iteration. These facts share a recurrence, but they do not serve the same role: the former supports finite coding, while the latter supplies the model of stable local self-application.
\end{remark}

\subsection{Mathematical limit}\label{subsec:mathematical-limit}

\begin{remark}
	The extremal property of $\Phi$ is used only in its number-theoretic sense. By \cite{hurwitz1891},
	\begin{equation}
		[\;1;1,1,\ldots\;]
	\end{equation}
	is the symbolic fixed point of the simplest continued-fraction constructor. Among irrational numbers, $\Phi$ is extremal for rational approximation in the Hurwitz/continued-fraction sense: its continued fraction has all partial quotients equal to $1$, making it the canonical worst-approximable irrational up to the sharp Hurwitz constant. This is a boundary model, not a theorem about proof search.
\end{remark}

\subsection{Carryless Pairing}\label{subsec:carryless-pairing}

\begin{definition}[Carryless pairing]
	Let $a,b\in\mathbb{N}$ have binary expansions
	\begin{equation}
		a=\sum_{i<k}\alpha_i2^i,\qquad b=\sum_{i<\ell}\beta_i2^i,\qquad \alpha_i,\beta_i\in\{0,1\}.
	\end{equation}
	The carryless pairing of $a$ and $b$ is
	\begin{equation}
		\pi_{\mathrm{CL}}(a,b)=\sum_i \alpha_i2^{2i}+\sum_i \beta_i2^{2i+1}.
	\end{equation}
	The first component occupies the even binary positions and the second component occupies the odd binary positions.
\end{definition}

\begin{remark}
	The construction is called carryless only in this restricted sense: the two components are placed on disjoint digit tracks, so no addition step can mix their bits. Decoding is obtained by projecting the binary expansion of $\pi_{\mathrm{CL}}(a,b)$ onto its even and odd positions. Thus $\pi_{\mathrm{CL}}$ is injective, and both projections are primitive recursive.
\end{remark}

\begin{remark}
	In the present paper, carryless pairing is only a technical substrate, discussed extensively elsewhere \citep{rosko2026carryless}. It supplies a finite code for ordered data such as line numbers, formula codes, and local proof-checking triples. It is not a new proof system and it carries no claim about the complexity of proof search.
\end{remark}

\subsection{Fibonacci Normal Forms}\label{subsec:fibonacci-normal-forms}

Let $(F_i)_{i\geq 1}$ be the shifted Fibonacci sequence
\begin{equation}
	F_1=1,\qquad F_2=2,\qquad F_{i+2}=F_{i+1}+F_i.
\end{equation}
A finite binary word $\epsilon=(\epsilon_1,\ldots,\epsilon_k)$ represents
\begin{equation}
	\langle\epsilon\rangle_F=\sum_{i=1}^{k}\epsilon_iF_i.
\end{equation}
The empty word represents $0$. The word is nonadjacent when
\begin{equation}
	\epsilon_i\epsilon_{i+1}=0\quad\text{for all }i<k.
\end{equation}
By Zeckendorf representation, every $n\in\mathbb{N}$ has a unique nonadjacent Fibonacci word representing it \citep{zeckendorf}.

\begin{proposition}[Arithmetical suitability]
	Fibonacci normalisation is primitive recursive. Equivalently, the relation saying that a finite binary word is the nonadjacent Fibonacci representation of $n$ is primitive recursive. Hence Fibonacci normalisation may be used as a canonical arithmetical coding layer for finite data.
\end{proposition}

\begin{proof}
	The Fibonacci function is primitive recursive. Given a finite word, its value is obtained by a bounded summation of Fibonacci terms, and nonadjacency is checked by a bounded scan of consecutive digit positions. The greedy construction of the Zeckendorf representation is a bounded primitive-recursive search for the largest Fibonacci term not exceeding the current remainder, followed by iteration on the strictly smaller remainder. Thus both evaluation and normal-form checking are primitive recursive.
\end{proof}

\begin{remark}
	The proof predicate below does not depend on Fibonacci normal forms. It works for any fixed primitive-recursive coding of formulas and finite sequences. Fibonacci normal forms are included only to show that the recurrence used in the boundary model can also be implemented as a canonical arithmetical coding discipline. Thus Fibonacci recurrence is suitable for arithmetization, but it is not an additional inference rule and not a source of proof-theoretic strength.
\end{remark}

\subsection{Pair Predicate}\label{subsec:pair-predicate}

\begin{definition}[Pair predicate]
	Let $\operatorname{even}(c)$ and $\operatorname{odd}(c)$ be the projections of the binary expansion of $c$ onto its even and odd positions. Define
	\begin{equation}
		\operatorname{Pair}(c,a,b)\quad\Longleftrightarrow\quad \operatorname{even}(c)=a\ \wedge\ \operatorname{odd}(c)=b.
	\end{equation}
	Equivalently, $\operatorname{Pair}(c,a,b)$ holds exactly when $c=\pi_{\mathrm{CL}}(a,b)$.
\end{definition}

\begin{remark}
	The projections are primitive recursive; in weak arithmetic, bounded representability depends on the chosen coding and available length and exponentiation functions. For the present external coding, the predicate gives a local check that a code really contains the ordered data it is supposed to contain.
\end{remark}

\begin{remark}
	$\operatorname{Pair}$ is syntactic. It verifies the shape of a finite code, not the truth of a formula and not the soundness of a proof system.
\end{remark}

\section{Proof Checking}\label{sec:proof-checking}

The proof-checking component is local and syntactic. A finite derivation is treated as a list of coded lines, and each line is accepted only when it is an axiom instance or follows from earlier lines by a bounded rule check. This establishes witnessed verification for a given proof object; it does not establish a global reflection principle for the system itself.

\subsection{Line Codes}\label{subsec:line-codes}

\begin{definition}[Line code]
	A proof line is coded as an ordered pair
	\begin{equation}
		\ell=\pi_{\mathrm{CL}}(t,\varphi),
	\end{equation}
	where $t$ is a tag and $\varphi$ is the code of the displayed formula. The tag records how the line is to be checked: as an axiom instance or as the result of a local inference rule. For an \textsc{MP} line, the tag also contains the cited indices $p,q<j$ of the two earlier lines to be compared.
\end{definition}

\begin{remark}
	Only the form of the code matters here. The formula code $\varphi$ may be supplied by any fixed primitive-recursive Gödel coding of formulas. The present construction does not depend on a special syntax for formulas, except that substitution and pattern matching for axiom schemata are decidable by bounded finite checks.
\end{remark}

\begin{remark}
	A derivation code is a finite sequence of line codes. Checking such a sequence means inspecting each line together with the earlier lines named by its tag. No line is accepted by appeal to the global soundness of the whole derivation.
\end{remark}

\subsection{Axiom Heads}\label{subsec:axiom-heads}

The displayed heads are for the implicational Hilbert fragment. They are not intended as axiom heads for full propositional logic \citep{hilbert28}.

\begin{definition}[Axiom schema pattern]
	An axiom head is a finite pattern identifying the outer syntactic form of an axiom schema. For the implicational fragment, the two standard heads are
	\begin{align}
		\mathbf{K}:&\quad A\to(B\to A),\\
		\mathbf{S}:&\quad (A\to(B\to C))\to((A\to B)\to(A\to C)).
	\end{align}
	A line with formula code $\varphi$ is accepted as an axiom instance when $\varphi$ matches one of these heads under a uniform substitution of formula codes for schematic letters.
\end{definition}

\begin{remark}
	This is a syntactic head check. It requires parsing the outer implication structure and verifying that repeated schematic letters receive the same formula code each time they occur. These are finite equality checks on subformula codes.
\end{remark}

\begin{remark}
	The axiom-head test by itself certifies only syntactic schema membership; semantic soundness enters only through the external local-soundness argument below.
\end{remark}

\subsection{MP Alignment}\label{subsec:mp-alignment}

\begin{definition}[MP alignment]
	A line $\ell_j$ is accepted by \textsc{Modus Ponens} alignment when its tag names earlier lines $\ell_p,\ell_q$ with $p,q<j$, and the checker parses the formula code on $\ell_q$ as
	\begin{equation}
		A\to B,
	\end{equation}
	compares the antecedent code $A$ with the formula code on $\ell_p$, and compares the consequent code $B$ with the formula code on $\ell_j$. The check is purely syntactic: it consists only of parsing the implication code and testing equality of the cited subformula codes.
\end{definition}

\begin{remark}
	This is the logical application pattern used in the resonance definition. It is not identical to proportional scaling. The proportional diagram keeps only the alignment idea: one object must match a distinguished part of another object, and a successor object is then determined.
\end{remark}

\begin{figure}[t]
	\makebox[\textwidth][c]{%
		\resizebox{0.95\textwidth}{!}{%
			\includegraphics[width=0.95\textwidth]{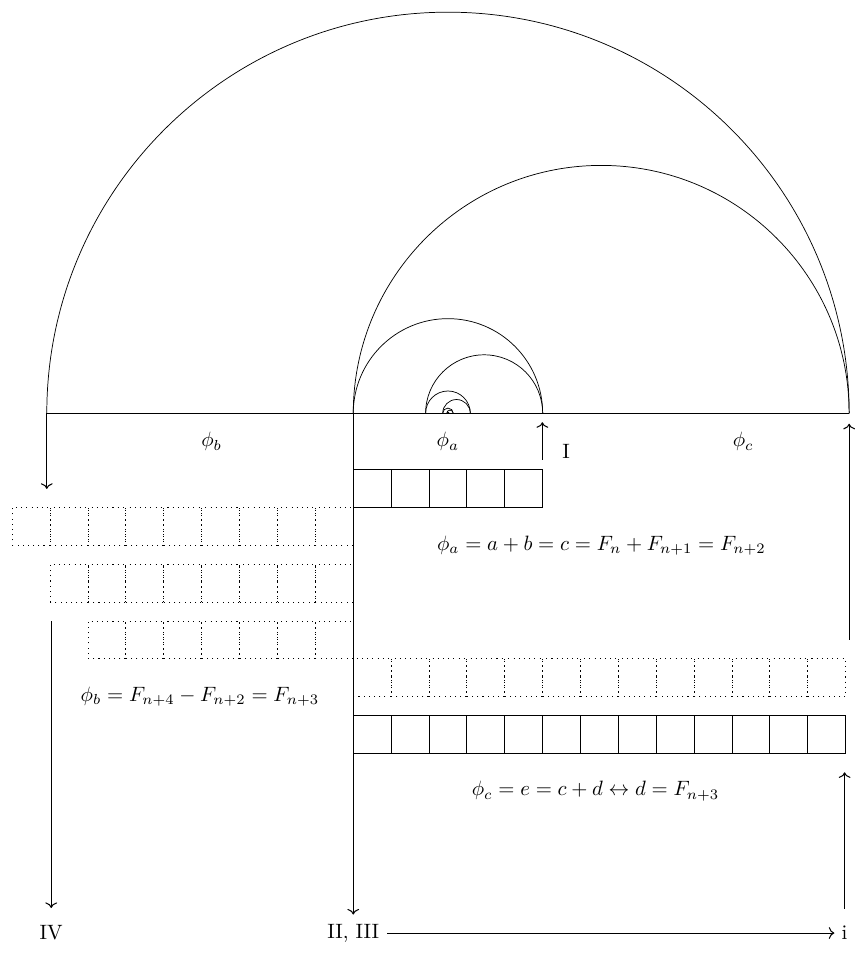}
		}%
	}
	\centering
	\caption{Geometric intuition as an analogy for local matching. The figure is an alignment diagram for the recurrence $F_n+F_{n+1}=F_{n+2}$. Three consecutive Fibonacci numbers $\{5,8,13\}$ appear as projections of arcs $\{\phi_a,\phi_b,\phi_c\}$ on nested circles scaled by $\Phi$ and $\bar{\Phi}$. The additional point $e=F_{n+4}$ illustrates the related check $F_{n+4}-F_{n+2}=F_{n+3}$: alignment preserves the recurrence pattern; deviation marks arithmetic failure.}
	\label{fig:iterant}
\end{figure}

\subsection{Bounded Witness}\label{subsec:bounded-witness}

\begin{definition}[Bounded witness]
	Let $D=(\ell_0,\ldots,\ell_{r-1})$ be a finite derivation code and let $\varphi$ be a formula code. Define
	\begin{equation}
		W(\varphi,D)
	\end{equation}
	to hold when the last line of $D$ has formula code $\varphi$, and every line of $D$ is accepted either by an axiom-head check or by an \textsc{MP} alignment using the earlier lines named in its tag.
\end{definition}

\begin{remark}
	The witness is bounded because every check is restricted to the finite object $D$. For line $j$, an \textsc{MP} tag supplies the cited pair $p,q<j$, and the verifier checks only those two earlier lines. The axiom check inspects only the finite formula code on the current line.
\end{remark}

\begin{remark}
	Thus $W(\varphi,D)$ says that $D$ is an explicit certificate for $\varphi$ relative to the chosen proof rules. It does not say that every true $\varphi$ has such a $D$, and it does not say that the system can certify its own global correctness.
\end{remark}

\section{Verification Scope}\label{sec:verification-scope}

The preceding definitions give a local verification procedure for finite proof objects. This section records the scope of that procedure. The point is not to derive a new logic, but to separate three levels: primitive-recursive checking, local soundness of accepted steps, and the absence of any internal global reflection principle.

\subsection{Primitive Recursion}\label{subsec:primitive-recursion}

\begin{proposition}
	The predicates used above to decode pairs, recognize axiom heads, test \textsc{MP} alignment, and verify $W(\varphi,D)$ are primitive recursive.
\end{proposition}

\begin{proof}
	Carryless decoding is projection onto even and odd binary positions. Axiom-head recognition is finite parsing plus equality of repeated subformula codes. \textsc{MP} alignment reads the cited indices from the tag and checks the two named earlier lines. The witness predicate $W(\varphi,D)$ is obtained by bounded iteration of these checks over the lines of $D$. Each operation is therefore a bounded computation on finite codes, hence primitive recursive.
\end{proof}

\begin{remark}
	Primitive recursiveness is an effectivity claim about the checker. It does not imply that short witnesses exist, nor that witnesses can be found efficiently.
\end{remark}

\subsection{Local Soundness}\label{subsec:local-soundness}

\begin{proposition}[Local soundness]
	Fix an intended semantics for the formula language. If each axiom head denotes a valid schema and \textsc{Modus Ponens} preserves validity under that semantics, then every line accepted by the checker is valid whenever the earlier cited lines are valid.
\end{proposition}

\begin{proof}
	The verification is by induction over the finite derivation code. An axiom line is valid by the assumed validity of the corresponding schema. An \textsc{MP} line cites earlier lines with codes $A$ and $A\to B$; by the induction hypothesis both are valid, and by soundness of \textsc{Modus Ponens} the consequent $B$ is valid. Thus each accepted line preserves validity relative to the earlier accepted lines.
\end{proof}

\begin{remark}
	This is local soundness only. It is a statement about the preservation of correctness along a given finite witness, not a claim that the system can prove its own soundness as a single internal theorem.
\end{remark}

\subsection{No Reflection}\label{subsec:no-reflection}

\begin{remark}
	The predicate $W(\varphi,D)$ is an external proof-checking predicate for particular finite derivations, relative to the chosen coding and rules. It does not yield an internal reflection principle such as $\mathrm{Prov}(\ulcorner\varphi\urcorner)\to\varphi$, nor does it yield global soundness or an internal theorem asserting the correctness of all accepted derivations. Such reflection principles are different proof-theoretic claims and are not obtained from primitive-recursive checking or local soundness; this is the standard proof-theoretic boundary associated with incompleteness and provability logic \citep{g31,lob55}.
\end{remark}

\section{Discussion}\label{sec:discussion}

The construction should be read as a boundary model. On one side, the reciprocal update $R(x)=1+1/x$ gives a disciplined form of local self-application: the same constructor is reapplied, the positive fixed point is unique, and finite approximations converge toward $\Phi$. On the other side, proof-theoretic reflection is not a local update. It asks a system to certify, uniformly and from within, the correctness of its own methods.

The proof-checking substrate developed above belongs to the local side of this boundary. A finite derivation can be decoded, inspected, and accepted line by line. Axiom heads and \textsc{MP} alignment are bounded syntactic checks, and a witness $D$ certifies a particular formula $\varphi$ relative to the rules. This is local certification, not global self-certification.

The role of $\Phi$ is therefore structural and contrastive. It does not govern logic, and Fibonacci recurrence does not provide a new proof system. Rather, $\Phi$ supplies a precise model of stable local self-application, while reflection marks the point where self-application is no longer justified by local witnessed checking alone. In this sense, the paper uses arithmetic not to reduce proof theory, but to make a boundary visible.

Finally, the Hurwitz extremality of $\Phi$ can be read as an obstruction template: where a system attempts to approximate a stable self-applicative limit by finite rational witnesses, $\Phi$ marks the case of maximal resistance to rational compression.

\clearpage
\vspace*{\fill}

\section*{References and Notes}

{\scriptsize
\bibliographystyle{plainnat}
\setlength{\bibsep}{0.1em}
\bibliography{refs}}

\subsection*{Note on this Version}

Earlier versions of this paper were unsatisfactory in presentation. They allowed the rhetoric of the construction to exceed what the formal components could responsibly support, especially around proof search, Fibonacci recurrence, and the role of $\Phi$ in proof theory. The present version has therefore been substantially overhauled. Its claims are restricted to a boundary model: local witnessed recursion is treated as stable self-application, global reflection as a distinct threshold.

\subsection*{Final Remarks}

The author welcomes criticism, proposed extensions, scholarly correspondence, and constructive dialogue. No conflicts of interest are declared. This research received no funding.

\begin{center}
	\scriptsize{
	\vspace{1em}
	Milan Rosko\\
	\vspace{1em}
	ORCID: \href{https://orcid.org/0009-0003-1363-7158}{\textsf{0009-0003-1363-7158}}\\[1ex]
	Email: \href{mailto:hi-at-milanrosko.com}{\textsf{hi-at-milanrosko.com}}\\[1ex]
	Email: \href{mailto:Q1012878@studium.fernuni-hagen.de}{\textsf{Q1012878@studium.fernuni-hagen.de}}\\
	\vspace{1em}
	Licensed under \ccby \\ \href{http://creativecommons.org/licenses/by/4.0/}{\scriptsize\textsf{creativecommons.org/licenses/by/4.0}}
	}
\end{center}
\vspace*{\fill}

\end{document}